\documentclass{article}%
\usepackage{graphicx}
\usepackage{amsmath}
\usepackage{amsfonts}
\usepackage{amssymb}%
\newtheorem{theorem}{Theorem}
{}

\newtheorem{lemma}{Lemma}
{}

\newenvironment{proof}[1][Proof]{\textbf{#1.} }{\ \rule{0.5em}{0.5em}}

\begin{document}

\title{On\ the Nonself-adjoint Sturm-Liouville Operator with Matrix Potential}
\author{O.A.Veliev\\{\small Depart. of Math, Fac.of Arts and Sci.,}\\{\small Dogus University, Ac\i badem, 34722,}\\{\small Kadik\"{o}y, \ Istanbul, Turkey.}\\\ {\small e-mail: oveliev@dogus.edu.tr}}
\date{}
\maketitle

\begin{abstract}
In this article we obtain the asymptotic formulas for the eigenvalues and
eigenfunctions of the nonself-adjoint operator generated by a system of
Sturm-Liouville equations with summable coefficients and the quasiperiodic
boundary conditions. Then using these asymptotic formulas, we find the
conditions on the potential for which the root funcions of this operator form
a Riesz basis.

\end{abstract}

\bigskip We consider the differential operator $L_{t}(Q(x))$ generated in the
space $L_{2}^{m}\left[  0,1\right]  $ of the vector functions by the
differential expression
\[
l(y)=-y^{^{\prime\prime}}(x)+Q\left(  x\right)  y(x)
\]
and the quasiperiodic conditions
\begin{equation}
y^{^{\prime}}\left(  1\right)  =e^{it}y^{^{\prime}}\left(  0\right)  ,\text{
}y\left(  1\right)  =e^{it}y\left(  0\right)
\end{equation}
for $t\in\lbrack0,2\pi),$ where $L_{2}^{m}\left[  0,1\right]  $ is the set of
the vector functions

$f\left(  x\right)  =\left(  f_{1}\left(  x\right)  ,f_{2}\left(  x\right)
,...,f_{m}\left(  x\right)  \right)  $ with $f_{k}\left(  x\right)  \in
L_{2}\left[  0,1\right]  $ for $k=1,2,...,m;$ and $Q\left(  x\right)  =\left(
b_{i,j}\left(  x\right)  \right)  $ is a $m\times m$ matrix with the
complex-valued summable entries $b_{i,j}\left(  x\right)  $. The norm
$\left\Vert .\right\Vert $ and inner product $(.,.)$ in $L_{2}^{m}\left[
0,1\right]  $ are defined by%
\[
\left\Vert f\right\Vert =\left(  \int\limits_{0}^{1}\left\vert f\left(
x\right)  \right\vert ^{2}dx\right)  ^{\frac{1}{2}},\text{ }(f\left(
x\right)  ,g\left(  x\right)  )=\int\limits_{0}^{1}\left\langle f\left(
x\right)  ,g\left(  x\right)  \right\rangle dx,
\]
where $\left\vert .\right\vert $ and $\left\langle .,.\right\rangle $ are the
norm and inner product in $\mathbb{C}^{m}.$ These boundary value problems play
a fundamental role in the spectral theory of the differential operator $L$
generated in the space $L_{2}^{m}\left(  -\infty,\infty\right)  $ by
expression $l(y)$ with periodic coefficient $Q\left(  x+1\right)  =Q\left(
x\right)  $, since the spectrum of the operator $L$ is the union of the
spectra of $L_{t}$ for $t\in\lbrack0,2\pi)$ ( see [1]).

Let us describe briefly the scheme of the paper.

First we show that ( see Theorem 1) it easily follows from the well-known
classical investigations [3] that the eigenvalues of the operator
$L_{t}\left(  Q\right)  $ consist of $m$ sequences
\begin{equation}
\{\lambda_{k,1}:k\in\mathbb{Z}\mathbb{\}},\text{ }\{\lambda_{k,2}%
:k\in\mathbb{Z}\mathbb{\}},...,\text{ }\{\lambda_{k,m}:k\in\mathbb{Z}%
\mathbb{\}}%
\end{equation}
lying in $O(|k|^{1-\frac{1}{m}})$ neighborhood of the eigenvalues $\left(
2k\pi+t\right)  ^{2}$ of the operator $L_{t}(0),$ where $L_{t}(Q)$ is denoted
by $L_{t}(0)$ if $Q(x)=0.$ Then we prove that the eigenvalues $\lambda_{k,l}$
of $L_{t}(q)$ lie in $O\left(  \frac{\ln|k|}{k}\right)  $ neighborhood of the
eigenvalues of the operator $L_{t}(C),$ where $C=\int\limits_{0}^{1}Q\left(
x\right)  dx.$ For this we consider the operator $L_{t}(Q)$ as perturbation of
$L_{t}(C)$\ by $Q(x)-C,$ that is, we take the operator$L_{t}(C)$ for an
unperturbed operator and the operator of multiplication by $Q(x)-C$ for a
perturbation. Therefore we analyze the eigenvalues and eigenfunction of
$L_{t}(C)$ and use the formulas (17)-(19) connecting the eigenvalues and
eigenfunctions of the operators $L_{t}(Q)$ and $L_{t}(C).$ Then we estimate
the terms in the connecting formulas ( see Lemma 2 and Lemma 3) by using Lemma
1. At last using the connecting formulas and lemmas 2, 3 we find asymptotic
formulas for eigenvalues and eigenfunctions of $L_{t}(Q)$ ( see Theorem 2).
Then using these asymptotic formulas, we prove that if the eigenvalues of the
matrix $C$ are simple then the root functions (eigenfunctions and associated
functions ) of the operator $L_{t}(Q)$ for $t\neq0,\pi$ form a Riesz basis.
The suggested method in this paper gives the possibility of obtaining the
asymptotic formulas of order $O(k^{-1}\ln|k|)$ for the eigenvalues
$\lambda_{k,j}$ and for the corresponding normalized eigenfunctions
$\Psi_{k,j}(x)$ of $L_{t}(Q)$ when the entries $b_{i,j}(x)$ of $Q(x)$ belong
to $L_{1}[0,1]$, that is, when there is not any condition about smoothness of
the\ coefficient $Q(x)$. Note that to obtain the asymptotic formulas of order
$O(\frac{1}{k})$ for the eigenvalues $\lambda_{k,j}$ of $L_{t}(Q)$ by using
the classical asymptotic expansions for the solutions of the matrix equation
$-Y^{^{\prime\prime}}+Q\left(  x\right)  Y=\lambda Y$ it is required \ that
$\ Q(x)$ be differentiable (see [3,4,6,7]). First we use the following theorem
which is easily obtained from the results of chapter 3 of [3]

\begin{theorem}
The boundary conditions (1) are regular and the eigenvalues of the operator
$L_{t}\left(  Q\right)  $ for $t\neq0,\pi$ consist of $m$ sequences (2)
satisfying
\begin{equation}
\lambda_{k,j}\left(  t\right)  =\left(  2\pi k+t\right)  ^{2}+O\left(
k^{1-\frac{1}{m}}\right)  \text{, }k=\pm N,\pm\left(  N+1\right)  ,...
\end{equation}
for $j=1,2,...,m$.
\end{theorem}

\begin{proof}
To prove that the conditions (1) are regular ( for definition see [3]) we need
to show that both the numbers $\theta_{-m}$ , $\theta_{m}$\ defined by the
equation
\begin{equation}
\theta_{-m}s^{-m}+\theta_{-m+1}s^{m-1}+...+\theta_{m}s^{m}=\det M(m),
\end{equation}
where $M(m)=\left[
\begin{array}
[c]{cc}%
\left(  e^{it}-s\right)  iI, & \left(  e^{it}-\frac{1}{s}\right)  \left(
-i\right)  I\\
\left(  e^{it}-s\right)  I, & \left(  e^{it}-\frac{1}{s}\right)  I
\end{array}
\right]  ,$

$I$ is $m\times m$ identity matrix, do not vanish. The right-hand side of (4)
is $2m\times2m$ determinant. First we prove that%
\begin{equation}
\det M(m)=\left(  -2ie^{it}s+2i+2ie^{2it}-2ie^{it}\frac{1}{s}\right)  ^{m}%
\end{equation}
Let us prove (5) by induction. This formula holds for $m=1$, since%
\begin{align}
\det M(1)  &  =\left\vert
\begin{array}
[c]{cc}%
\left(  e^{it}-s\right)  i, & \left(  e^{it}-\frac{1}{s}\right)  \left(
-i\right) \\
\left(  e^{it}-s\right)  , & \left(  e^{it}-\frac{1}{s}\right)
\end{array}
\right\vert \\
&  =-2ie^{it}s+2i+2ie^{2it}-2ie^{it}\frac{1}{s}.\nonumber
\end{align}
Assume that (5) holds for $m=p-1$, that is,%
\begin{equation}
\det M(p-1)=\left(  -2ie^{it}s+2i+2ie^{2it}-2ie^{it}\frac{1}{s}\right)
^{p-1}.
\end{equation}
Then using the obvious equality
\[
\det M(p)=\left\vert
\begin{tabular}
[c]{ll}%
$M_{1}$ & $0$\\
$0$ & $M_{p-1}$%
\end{tabular}
\ \ \right\vert
\]
we see that $\det M(p)=(\det M(1))\det M(p-1).$ Hence (6) and (7) imply (5).
It follows from\ (4) and (5) that%
\[
\theta_{m}=\left(  -2ie^{it}\right)  ^{m}\neq0\text{, }\theta_{-m}=\left(
-2ie^{it}\right)  ^{-m}\neq0,\text{ }\forall t\in\lbrack0,2\pi)
\]
which means that the boundary conditions (1) are regular. By (5) the numbers
$s_{1}=e^{it}$ and $s_{2}=e^{-it}$ for $t\neq0,\pi$ are roots of the equation%
\begin{equation}
\det M(m)=0
\end{equation}
with multiplicity $m$, since $e^{it}\neq e^{-it}$. Therefore it follows from
theorem 2 of chapter 3 of [3] that to each multiple zero of (8) ,that is, to
$s_{1}$and $s_{2}$ correspond $m$ sequences

$\{\lambda_{k,1}^{(1)}:k=N,N+1,...\mathbb{\}}$, $\{\lambda_{k,2}%
^{(1)}:k=N,N+1,...\mathbb{\}}$,...,$\{\lambda_{k,m}^{(1)}%
:k=N,N+1,...\mathbb{\}}$

and

$\{\lambda_{k,1}^{(2)}:k=N,N+1,...\mathbb{\}}$, $\{\lambda_{k,2}%
^{(2)}:k=N,N+1,...\mathbb{\}}$,...,$\{\lambda_{k,m}^{(2)}%
:k=N,N+1,...\mathbb{\}}$ respectively, satisfying%
\[
\lambda_{k,j}^{(1)}=\left(  2k\pi+t\right)  ^{2}+O\left(  k^{1-\frac{1}{m}%
}\right)  \text{, }k=N,N+1,...,
\]%
\[
\lambda_{k,j}^{(2)}=\left(  2k\pi-t\right)  ^{2}+O\left(  k^{1-\frac{1}{m}%
}\right)  \text{, }k=N,N+1,...
\]
for $j=1,2,...,m.$ Clearly, these $2$ formulas can be written in the form (3)
by taking $k=\pm N,\pm\left(  N+1\right)  ,...$ .
\end{proof}

In forthcoming relations by $N$ we denote a big integer, that is, $N\gg1$, and
by $c_{m}$, for $m=1,2,...$, the positive constants, independent on $N,$ whose
exact values are inessential. The formula (3) shows that the eigenvalue
$\lambda_{k,j}$ of the operator $L_{t}(Q)$ for $t\neq0,\pi$ is close to the
eigenvalue $\left(  2k\pi+t\right)  ^{2}$ \ of the operator $L_{t}(0)$ and far
from the other eigenvalues $\left(  2n\pi+t\right)  ^{2}$ \ ( $n\neq k$) of
$L_{t}(0).$ Namely, for $t\neq0,\pi$ we have the inequalities
\[
|\lambda_{k,j}-\left(  2k\pi+t\right)  ^{2}|<c_{1}|k|^{1-\frac{1}{m}},\text{
}|\lambda_{k,j}-\left(  2(-k)\pi+t\right)  ^{2}|>c_{2}|k|,
\]%
\begin{equation}
|\lambda_{k,j}-\left(  2\pi n+t\right)  ^{2}|>c_{2}%
(||k|-|n||)(|k|+|n|),\forall n\neq\pm k.
\end{equation}
These inequalities imply the following purposive relations%
\begin{equation}
\sum_{n:n>d}\frac{1}{\left\vert \lambda_{k,j}-(2\pi n+t)^{2}\right\vert
}<\frac{c_{3}}{d},\text{ }\forall d>2\mid k\mid,
\end{equation}%
\begin{equation}
\sum_{n:n\neq k}\frac{1}{\left\vert \lambda_{k,j}-(2\pi n+t)^{2}\right\vert
}=O((\frac{\ln|k|}{k}),
\end{equation}%
\begin{equation}
\sum_{n:n\neq k}\frac{1}{\left\vert \lambda_{k,j}-(2\pi n+t)^{2}\right\vert
^{2}}=O(\frac{1}{k^{2}})
\end{equation}
for $|k|\geq N$ and $t\neq0,\pi$. Clearly,

$\varphi_{n,1}=\left(
\begin{array}
[c]{c}%
e^{i\left(  2\pi n+t\right)  x}\\
0\\
\vdots\\
0
\end{array}
\right)  ,\varphi_{n,2}=\left(
\begin{array}
[c]{c}%
0\\
e^{i\left(  2\pi n+t\right)  x}\\
\vdots\\
0
\end{array}
\right)  ,...,\varphi_{n,m}=\left(
\begin{array}
[c]{c}%
0\\
\vdots\\
0\\
e^{i\left(  2\pi n+t\right)  x}%
\end{array}
\right)  $ are the eigenfunctions of the operator $L_{t}(0)$ corresponding to
the eigenvalue $\left(  2\pi n+t\right)  ^{2}$. The multiplicity of the
eigenvalue $\left(  2\pi n+t\right)  ^{2}$ is $m$ and the corresponding
eigenspace is $E_{n}=span\left\{  \varphi_{n,1},\varphi_{n,2},...,\varphi
_{n,m}\right\}  $. The adjoint operator to $L_{t}\left(  Q\right)  $ is
$L_{t}\left(  Q^{\ast}\right)  ,$ where $Q^{\ast}\left(  x\right)  $ is the
adjoint matrix to $Q\left(  x\right)  $, that is,%

\[
Q^{\ast}\left(  x\right)  =\left[
\begin{array}
[c]{cccc}%
\overline{b_{1,1}}\left(  x\right)  & \overline{b_{2,1}}\left(  x\right)  &
... & \overline{b_{m,1}}\left(  x\right) \\
\overline{b_{1,2}}\left(  x\right)  & \overline{b_{2,2}}\left(  x\right)  &
... & \overline{b_{m,2}}\left(  x\right) \\
\vdots & \vdots & \vdots & \vdots\\
\overline{b_{1,m}}\left(  x\right)  & \overline{b_{2,m}}\left(  x\right)  &
... & \overline{b_{m,m}}\left(  x\right)
\end{array}
\right]
\]
Since the boundary conditions (1) are self-adjoint the operator $L_{t}\left(
Q\right)  $ is selfadjoint if $Q(x)$ is a symmetric matrix. Hence $L_{t}%
^{\ast}(0)=L_{t}(0).$

Now to analyze the operator $L_{t}(C)$ we introduce the following notations.
Suppose the matrix $C$ has $p$ distinct eigenvalues $\mu_{1},\mu_{2}%
,...,\mu_{p}$ with multiplicities $m_{1},m_{2},...,m_{p}$ respectively (
$m_{1}+m_{2}+...+m_{p}=m)$. Let $u_{j,1},$ $u_{j,2},...u_{j,s_{j}}$ be the
eigenvectors corresponding to the eigenvalue $\mu_{j}$ and $u_{j,s,1},$
$u_{j,s,2},...u_{j,s,r_{j,s}-1}$ be associated vectors belonging to the
eigenvector $u_{j,s}.$ Note that $r_{j,s}$ is called the multiplicity of the
eigenfunction $u_{j,s}$ and $r_{j,1}+r_{j,2}+...+r_{j,s_{j}}=m_{j}.$ The
number $r_{j}$ defined by
\begin{equation}
r_{j}=\max_{s}r_{j,s}.
\end{equation}
is a maximum multiplicity of the eigenfunctions corresponding to the
eigenvalue $\mu_{j}.$ It is not hard to verify that $u_{j,s}e^{i\left(  2\pi
k+t\right)  x}$ for $s=1,2,...,$ are the eigenfunction of $L(C)$ corresponding
to the eigenvalue $\mu_{k,j}=\left(  2\pi k+t\right)  ^{2}+\mu_{j}$ and
$u_{j,s,r}e^{i\left(  2\pi k+t\right)  x}$ for $r=1,2,...,$ are the associated
function of $L(C)$ belonging to $u_{j,s}e^{i\left(  2\pi k+t\right)  x}$ .

In subsequent relations, for simplicity of notation, the eigenvalues of $C,$
counted with multiplicity, is indexed as $\mu_{1},\mu_{2},...,\mu_{m}$. Any
normalized eigenvector corresponding to the eigenvalue $\mu_{j}$ is denoted by
$v_{j}$ and the associated vectors belonging to the eigenvector $v_{j}$ are
denoted by $v_{j,s}$ for $s=1,2,...$ In this notation the eigenvalues,
eigenfunction, associated function of $L_{t}(C)$ are

$\mu_{k,j}=\left(  2\pi k+t\right)  ^{2}+\mu_{j},$ $\Phi_{k,j}(x)=v_{j}%
e^{i\left(  2\pi k+t\right)  x},$ $\Phi_{k,j,s}(x)=v_{j,s}e^{i\left(  2\pi
k+t\right)  x}$ respectively. Similarly, the eigenvalues, eigenfunction and
associated function of $L_{t}^{\ast}(C)$ are $\overline{\mu_{k,j}},$
$\Phi_{k,j}^{\ast}(x)=v_{j}^{\ast}e^{i\left(  2\pi k+t\right)  x},$ and
$\Phi_{k,j,s}^{\ast}(x)=v_{j,s}^{\ast}e^{i\left(  2\pi k+t\right)  x}$ , where
$v_{j}^{\ast},$ $v_{j,s}^{\ast}$ are the eigenvector and associated vector of
\ $C^{\ast}$ corresponding to $\overline{\mu_{j}}.$ By definition we have
\begin{align}
(L^{\ast}(C)-\overline{\mu_{k,j}})\Phi_{k,j}^{\ast}(x)  &  =0,\\
(L^{\ast}(C)-\overline{\mu_{k,j}})\Phi_{k,j,s}^{\ast}(x)  &  =\Phi
_{k,j,s-1}^{\ast}(x),
\end{align}
where $\Phi_{k,j,0}^{\ast}(x)=\Phi_{k,j}^{\ast}(x)$. Multiplying both sides
of
\begin{equation}
L(Q\left(  x\right)  )\Psi_{k,j}(x)=\lambda_{k,j}\Psi_{k,j}(x)
\end{equation}
by $\Phi_{k,j}^{\ast}(x),$ using $L(Q)=L(C)+(Q(x)-C)$ and (14), we get%
\begin{equation}
(\lambda_{k,j}-\mu_{k,j})(\Psi_{k,j}(x),\Phi_{k,j}^{\ast}(x))=((Q(x)-C)\Psi
_{k,j}(x),\Phi_{k,j}^{\ast}(x)).
\end{equation}
Now multiplying (16) by $\Phi_{k,j,1}^{\ast}(x)$ and using (15), (17), we get%
\[
(\lambda_{k,j}-\mu_{k,j})^{2}(\Psi_{k,j}(x),\Phi_{k,j,1}^{\ast}%
(x))=((Q(x)-C)\Psi_{k,j}(x),\Phi_{k,j}^{\ast}(x))+
\]%
\begin{equation}
(\lambda_{k,j}-\mu_{k,j})((Q(x)-C)\Psi_{k,j}(x),\Phi_{k,j,1}^{\ast}(x)).
\end{equation}
In this way one can deduce the formulas
\begin{equation}
(\lambda_{k,j}-\mu_{k,j})^{s+1}(\Psi_{k,j}(x),\Phi_{k,j,s}^{\ast}(x))=
\end{equation}%
\[
\sum_{p=0}^{s}(\lambda_{k,j}-\mu_{k,j})^{p}((Q(x)-C)\Psi_{k,j}(x),\Phi
_{k,j,p}^{\ast}(x)).
\]
To estimate the terms $((Q(x)-C)\Psi_{k,j}(x),\Phi_{k,j,p}^{\ast}(x)),$
$(\Psi_{k,j}(x),\Phi_{k,j,s}^{\ast}(x))$ of the formula (19)\ we use
(10)-(12), the following lemma, and the formulas
\begin{equation}
\left(  \lambda_{k,j}-\left(  2\pi n+t\right)  ^{2}\right)  \left(  \Psi
_{k,j}(x),\varphi_{n,s}(x)\right)  =\left(  \Psi_{k,j}(x),Q^{\ast}\left(
x\right)  \varphi_{n,s}(x)\right)  ,
\end{equation}

\begin{equation}
\left(  \Psi_{k,j}(x),\varphi_{n,s}(x)\right)  =\frac{\left(  \Psi
_{k,j}(x),Q^{\ast}\left(  x\right)  \varphi_{n,s}(x)\right)  }{\lambda
_{k,j}-\left(  2\pi n+t\right)  ^{2}},\text{ }\forall n\neq k
\end{equation}
which can be obtained from(16) by multiplying both sides by $\varphi_{n,s}(x)$
and using $L_{t}\left(  0\right)  \varphi_{n,s}(x)=\left(  2\pi n+t\right)
^{2}\varphi_{n,s}(x),$ $L_{t}^{\ast}\left(  0\right)  =L_{t}\left(  0\right)
$.

\begin{lemma}
If $t\neq0,\pi,$ then
\begin{equation}
\left(  \Psi_{k,j}(x),Q^{\ast}\left(  x\right)  \varphi_{n,s}(x)\right)
=\sum\limits_{\substack{q=1,2,...m\\p=-\infty,...,\infty}}b_{s,q,n-p}%
(\Psi_{k,j}(x),\varphi_{p,q}(x)),
\end{equation}%
\begin{equation}
\text{ }\left\vert \left(  \Psi_{k,j}(x),Q^{\ast}\left(  x\right)
\varphi_{n,s}(x)\right)  \right\vert <c_{4}\text{ }%
\end{equation}
for $n\in Z$ $;$ $\mid k\mid\geq N$ $;$ $s,j=1,2,...,m$, where $\ b_{s,q,n-p}%
=\int\limits_{0}^{1}b_{s,q}\left(  x\right)  e^{2\pi i\left(  p-n\right)
x}dx.$
\end{lemma}

\begin{proof}
Since $Q(x)\Psi_{k,j}(x)\in L_{1}^{m}[0,1]$ we have
\[
\lim_{n\rightarrow\infty}\left(  Q(x)\Psi_{k,j}(x),\varphi_{n,s}(x)\right)
=0,\text{ }\forall s=1,2,...,m.
\]
Therefore there exists a positive constant $C(k,j)$ and indices $n_{0},$
$j_{0}$ satisfying
\begin{equation}
\max_{\substack{n\in\mathbb{Z},\\s=1,2,...,m}}\left\vert (\Psi_{k,j},Q^{\ast
}\left(  x\right)  \varphi_{n,s})\right\vert =\left\vert (\Psi_{k,j},Q^{\ast
}\left(  x\right)  \varphi_{n_{0},j_{0}})\right\vert =C(k,j)
\end{equation}
Using this, (21), and then (10), we obtain
\begin{equation}
\mid\left(  \Psi_{k,j}(x),\varphi_{n,s}(x)\right)  \mid\leq\frac{C(k,j)}%
{\mid\lambda_{k,j}-\left(  2\pi n+t\right)  ^{2}\mid}%
\end{equation}%
\[
\sum_{n:n>d}\mid\left(  \Psi_{k,j}(x),\varphi_{n,s}(x)\right)  \mid
<\frac{c_{3}C(k,l)}{d},
\]
where $d>2|k|.$ This implies that the decomposition of $\Psi_{k,j}(x)$ by the
orthonormal basis $\{\varphi_{n,s}(x)$:$n\in\mathbb{Z},$ $s=1,2,...,m\}$ is of
the form
\begin{equation}
\Psi_{k,j}(x)=\sum_{\substack{p:|p|\leq d,\\q=1,2,...,m}}\left(  \Psi
_{k,j}(x),\varphi_{p,q}(x)\right)  \varphi_{p,q}(x)+g_{d}(x),
\end{equation}
where $\sup_{x\in\lbrack0,1]}|g_{d}(x)|<\frac{c_{3}C(k,l)}{d}$. Putting these
in $\left(  \Psi_{k,j}(x),Q^{\ast}\left(  x\right)  \varphi_{n,s}(x)\right)  $
and tending $d$ to $\infty$, we obtain (22). \ 

Now we prove (23). Using (26) in $(\Psi_{k,j}(x),Q^{\ast}\left(  x\right)
\varphi_{n_{0},j_{0}}(x))$, tending $d$ to $\infty,$ isolating the terms with
multiplicands $(\Psi_{k,j}(x),\varphi_{k,i}(x))$ for $i=1,2,...,m$ we obtain
\begin{align}
\left(  \Psi_{k,j}(x),Q^{\ast}\left(  x\right)  \varphi_{n_{0},j_{0}%
}(x)\right)   &  =\sum\limits_{i=1,2,...,m}b_{j_{0},i,n_{0}-k}(\Psi
_{k,j},\varphi_{k,i})+\nonumber\\
&  \sum\limits_{\substack{n:n\neq k\\i=1,2,...,m}}b_{j_{0},i,n_{0}-n}\left(
\Psi_{k,j},\varphi_{n,i}\right)
\end{align}
Since
\begin{equation}
\mid b_{j,i,s}\mid\leq\max_{p,q=1,2,...,m}\int\limits_{0}^{1}\mid
b_{p,q}\left(  x\right)  \mid dx<c_{5},\forall j,i,s
\end{equation}
it follows from (25) and (11) \ that
\[
\sum\limits_{\substack{n:n\neq k\\i=1,2,...,m}}b_{j_{0},i,n_{0}-n}\left(
\Psi_{k,j},\varphi_{n,i}\right)  =O(\frac{\ln|k|}{k}C(k,j)).
\]
Therefore taking into account that the absolute value of the first summation
in the right-hand side of (27) is not greater than $mc_{5}$ ( see (28)), we
conclude that

$\mid C(k,l)\mid<2mc_{5}$ ( see (24)) which means that (23) holds
\end{proof}

\begin{lemma}
If $t\neq0,\pi$ then
\begin{equation}
\left(  \Psi_{k,j}(x),(Q^{\ast}\left(  x\right)  -C^{\ast})\Phi_{k,i,p}^{\ast
}(x)\right)  =O(\frac{\ln|k|}{k})
\end{equation}
for $i=1,2,...,m$ and $p=0,1,2,....$
\end{lemma}

\begin{proof}
\bigskip Since $\Phi_{k,n,p}^{\ast}(x)\equiv v_{n,p}^{\ast}e^{i\left(  2\pi
k+t\right)  x}$ it is enough to prove that%
\begin{equation}
\left(  \Psi_{k,j}(x),(Q^{\ast}\left(  x\right)  -C^{\ast})\varphi
_{k,s}(x)\right)  =O((\frac{\ln|k|}{k}).
\end{equation}
for $s=1,2,...,m.$ Using the obvious relation
\begin{equation}
\left(  \Psi_{k,j}(x),C^{\ast}\varphi_{k,s}(x)\right)  =\sum
\limits_{i=1,2,...,m}b_{s,i,0}(\Psi_{k,j}(x),\varphi_{k,i}(x))\nonumber
\end{equation}
and (22), we see that
\begin{equation}
\left(  \Psi_{k,j},(Q^{\ast}\left(  x\right)  -C^{\ast})\varphi_{k,s}%
(x)\right)  =\sum\limits_{\substack{n:n\neq k\\i=1,2,...,m}}b_{s,i,k-n}\left(
\Psi_{k,j}(x),\varphi_{n,i}(x)\right)  .
\end{equation}
On the other hand, it follows from (21), (23) that
\begin{equation}
\mid\left(  \Psi_{k,j}(x),\varphi_{n,i}(x)\right)  \mid\leq\frac{c_{4}}%
{\mid\lambda_{k,j}-\left(  2\pi n+t\right)  ^{2}\mid}%
\end{equation}
for $n\neq k$ and $i=1,2,...,m.$ Therefore using this and (28), (11) we see
that the right-hand side of (31) is $O((\frac{\ln|k|}{k}).$ Hence (30) is proved
\end{proof}

\begin{lemma}
For each eigenfunction $\Psi_{k,j}(x)$ of $L_{t}(Q)$ for $|k|\geq N$ there
exists a root function $\Phi_{k,i,s}^{\ast}(x)$ of $L_{t}(C^{\ast})$
satisfying
\[
\left\vert \left(  \Psi_{k,j}(x),\Phi_{k,i,s}^{\ast}(x)\right)  \right\vert
>c_{6}.
\]

\end{lemma}

\begin{proof}
Since $\left\{  \varphi_{n,i}(x):i=1,2,...,m\text{, }n\in Z\right\}  $ is
orthonormal basis in $L_{2}^{m}\left[  0,1\right]  $ we have%
\[
\Psi_{k,j}(x)=\sum\limits_{i=1,2,...,m}\text{ }(\sum_{n\in Z}\left(
\Psi_{k,j}(x),\varphi_{n,i}(x)\right)  \varphi_{n,i}(x)),
\]%
\begin{equation}
\left\Vert \Psi_{k,j}\right\Vert ^{2}=\sum\limits_{i=1,2,...,m}\left\vert
\left(  \Psi_{k,j},\varphi_{k,i}\right)  \right\vert ^{2}+\sum
\limits_{i=1,2,...,m}\text{ }(\sum_{n:n\neq k}\left\vert \left(  \Psi
_{k,j},\varphi_{n,i}\right)  \right\vert ^{2})
\end{equation}
On the other hand it follows from (32) and (12) that
\begin{equation}
\sum\limits_{i=1,2,...,m}\text{ }(\sum_{n:n\neq k}\left\vert \left(
\Psi_{k,j},\varphi_{n,i}\right)  \right\vert ^{2})=O(\frac{1}{k^{2}})
\end{equation}
Since $\left\Vert \Psi_{k,j}(x)\right\Vert =1$, from (33) we obtain%
\begin{equation}
\sum\limits_{i=1,2,...,m}\left\vert \left(  \Psi_{k,j}(x),\varphi
_{k,i}(x)\right)  \right\vert ^{2}=1+O\left(  \frac{1}{k^{2}}\right)
\end{equation}
Thus the norm of the projection of $\Psi_{k,j}(x)$ on the subspace

$E_{k}=span\left\{  \varphi_{k,1},\varphi_{k,2},...,\varphi_{k,m}\right\}  $
is $1+O\left(  \frac{1}{k^{2}}\right)  .$ Therefore taking into account that
$\Phi_{k,j}^{\ast}(x)=v_{j}^{\ast}e^{i\left(  2\pi k+t\right)  x}$ ,
$\Phi_{k,j,s}^{\ast}(x)=v_{j,s}^{\ast}e^{i\left(  2\pi k+t\right)  x}$ and
$v_{j},v_{j,s}^{\ast}$ for $j=1,2,...;$ $s=1,2,...,$ are the system of
eigenvectors and associated vectors of \ $C^{\ast},$ that is, form a basis in
$\mathbb{C}^{m},$ we get the prove of the lemma
\end{proof}

\begin{theorem}
\bigskip(Main). If $t\neq0,\pi$ then:

$(a)$ All big eigenvalues of $L_{t}\left(  Q\left(  x\right)  \right)  $ lie
in $O((\frac{\ln|k|}{k})^{\frac{1}{r_{j}}})$ neighborhood of the eigenvalues
$\mu_{k,j}=\left(  2\pi k+t\right)  ^{2}+\mu_{j}$ for $k\in\mathbb{Z}$,
$j=1,2,...,m$ of $L_{t}(C),$ where $r_{j}$ is defined in (13) and $r_{j}=1$ if
the matrix $C$ has no associated function corresponding to the eigenvalue
$\mu_{j}.$

$(b)$ Let $\mu_{j}$ be a simple eigenvalue of the matrix $C$ and
$\lambda_{k,j}$ be an eigenvalue of $L_{t}\left(  Q\left(  x\right)  \right)
$ lying in $\frac{1}{2}a_{j}$ neighborhood of $\mu_{k,j}=\left(  2\pi
k+t\right)  ^{2}+\mu_{j},$ where

$a_{j}=\min_{i\neq j}\mid\mu_{j}-\mu_{i}\mid.$ Then $\lambda_{k,j}$ is the
simple eigenvalue of $L_{t}\left(  Q\left(  x\right)  \right)  .$ Moreover
$\lambda_{k,j}$ and the corresponding eigenfunction $\Psi_{k,j}(x)$ satisfy
\begin{equation}
\lambda_{k,j}(t)=\left(  2\pi k+t\right)  ^{2}+\mu_{j}+O(\frac{\ln|k|}{k}).
\end{equation}%
\begin{equation}
\Psi_{k,j}(x)=v_{j}e^{i\left(  2\pi k+t\right)  x}+O(\frac{\ln|k|}{k}),
\end{equation}
where $v_{j}$ is the eigenvector of $C$ corresponding to the eigenvalue
$\mu_{j}.$

$(c)$ Suppose that the all eigenvalues $\mu_{1},\mu_{2},...,\mu_{m}$ of $C$
are simple. Then there exist a number $N$ such that all eigenvalues
$\lambda_{k,1},\lambda_{k,2}$,...,$\lambda_{k,m}$ of $L_{t}\left(  Q\left(
x\right)  \right)  $ for $\mid k\mid\geq N$ are simple and satisfy the
asymptotic formulas (36). The eigenfunction $\Psi_{k,j}(x)$ of $L(Q)$
satisfies (37). The root functions of $L_{t}\left(  Q\left(  x\right)
\right)  $ form a \ Riesz basis in $L_{2}^{m}(0,1).$
\end{theorem}

\begin{proof}
In (19) replacing $\mu_{k,j}$ and $\Phi_{k,j,n}^{\ast}(x)$ with $\mu_{k,i}$
and $\Phi_{k,i,n}^{\ast}(x)$ for $n=0,1,...,s,$ we get%
\begin{align}
&  (\lambda_{k,j}-\mu_{k,i})^{s+1}(\Psi_{k,j}(x),\Phi_{k,i,s}^{\ast}(x))\\
&  =\sum_{p=0}^{s}(\lambda_{k,j}-\mu_{k,i})^{p}((Q(x)-C)\Psi_{k,j}%
(x),\Phi_{k,i,p}^{\ast}(x))\nonumber
\end{align}
This formula, (29), and Lemma 3 imply that
\[
(\lambda_{k,j}-\mu_{k,i})^{s+1}=\sum_{p=0}^{s}(\lambda_{k,j}-\mu_{k,i}%
)^{p}O(\frac{\ln|k|}{k}),
\]
where $s+1\leq r_{i}$ (see (13)), from which we obtain the proof of \ $(a).$

$(b)$ Since $\mid\lambda_{k,j}-\mu_{k,i}\mid>\frac{1}{2}a_{j}$ for $i\neq j,$
using (29) and (38), we get%
\[
(\Psi_{k,j}(x),\Phi_{k,i,p}^{\ast}(x))=O\left(  \frac{\ln|k|}{k}\right)
\]
for $i\neq j$ and $p=1,2,....$ This and (34) imply that (37) holds for any
normalized eigenfunction corresponding to $\lambda_{k,j},$ since

$span\left\{  \varphi_{k,1},\varphi_{k,2},...,\varphi_{k,m}\right\}
=span\{\Phi_{k,i,p}^{\ast}(x):i=1,2,...;p=0,1,...\}$ and $\mu_{j}$ is a simple
eigenvalue of $C.$ Using this let us prove that $\lambda_{k,j}$ is a simple
eigenvalue. Suppose to the contrary that $\lambda_{k,j}$ is a multiple
eigenvalue. If there are two linearly independent eigenfunctions corresponding
to $\lambda_{k,j}$ then one can find two orthogonal eigenfunctions satisfying
(37) which is impossible. Hence there exists unique eigenfunction $\Psi
_{k,j}(x)$ corresponding to $\lambda_{k,j}.$ If there exists associated
function $\Psi_{k,j,1}(x)$ belonging to the eigenfunction $\Psi_{k,j}(x)$
then
\[
(L(Q)-\lambda_{k,j})\Psi_{k,j,1}(x)=\Psi_{k,j}(x).
\]
Multiplying both sides of this equality by $\Psi_{k,j}^{\ast}(x),$ where
$\Psi_{k,j}^{\ast}(x)$ is the eigenfunction of $L^{\ast}(Q)$ corresponding to
the eigenvalue $\overline{\lambda_{k,j}},$ we obtain%
\begin{equation}
(\Psi_{k,j},\Psi_{k,j}^{\ast})=(\Psi_{k,i,1}(x),(L^{\ast}(Q)-\overline
{\lambda_{k,j}}))\Psi_{k,j}^{\ast}(x))=0.
\end{equation}
Since the proved statements are also applicable for adjoint operator

$L^{\ast}(Q)=L(Q^{\ast}),$ it follows from (37) that
\begin{equation}
\Psi_{k,j}^{\ast}(x)=v_{j}^{\ast}e^{i\left(  2\pi k+t\right)  x}+O\left(
\frac{\ln|k|}{k}\right)  ,
\end{equation}
where $v_{j}^{\ast}$ is the eigenvector of $C^{\ast}$ corresponding to
$\overline{\mu_{j}}$. This formula with (37) and the obvious relation
$(v_{j},v_{j}^{\ast})\neq0$ contradicts (39). Thus $\lambda_{k,j}$ is a simple
eigenvalue and (37) hold. The formula (36) follows from $(a),$ since
$r_{j}=1.$

$(c)$ In $(a)$ we proved that all large eigenvalues of $L_{t}\left(  Q\left(
x\right)  \right)  $ lie in $\frac{c_{7}\ln|k|}{k}$ neighborhood $\Delta
_{k,l}=\{z\in\mathbb{C}$: $\mid z-\mu_{k,j}\mid<\frac{c_{7}\ln|k|}{k}\}$ of
$\mu_{k,j}=\left(  2\pi k+t\right)  ^{2}+\mu_{j}$ for $\mid k\mid\geq N$,
$j=1,2,...,m.$ Clearly, the circles $\Delta_{k,l}$ for $j=1,2,...,m$ and $\mid
k\mid>N$ are pairwise disjoint. Let us we prove that each of these circles
does not contain more than one eigenvalues of \ $L_{t}\left(  Q\left(
x\right)  \right)  .$ Suppose to the contrary that two different eigenvalues
$\Lambda_{1},\Lambda_{2}$ lie in $\Delta_{k,l}.$ Then by $(b)$ these
eigenvalues are simple and the corresponding eigenfunctions $\Psi_{1},\Psi
_{2}$ satisfy
\[
\Psi_{p}(x)=v_{j}e^{i\left(  2\pi k+t\right)  x}+O\left(  \frac{\ln|k|}%
{k}\right)
\]
for $p=1,2.$ Similarly, the eigenfunctions $\Psi_{1}^{\ast},\Psi_{2}^{\ast}$
of $L_{t}^{\ast}\left(  Q\left(  x\right)  \right)  $ corresponding to the
eigenvalues $\overline{\Lambda_{1}},\overline{\Lambda_{2}}$ satisfy
\[
\Psi_{p}^{\ast}(x)=v_{j}^{\ast}e^{i\left(  2\pi k+t\right)  x}+O\left(
\frac{\ln|k|}{k}\right)  .
\]
for $p=1,2.$ Since $\Lambda_{1}\neq\Lambda_{2}$ we have $0=(\Psi_{1}%
(x),\Psi_{2}^{\ast}(x))=(v_{j},v_{j}^{\ast})+O\left(  \frac{\ln|k|}{k}\right)
,$ which is impossible. Thus by Theorem 1 and $(a)$ the pairwise disjoint
circles $\Delta_{k,1},$ $\Delta_{k,2},...,$ $\Delta_{k,m}$ , where $|k|\geq
N,$ contain $m$ eigenvalues of $L_{t}(Q)$. On the other hand we just proved
that each of these circles does not contain more than one eigenvalues. Hence
for $j=1,2,...,m$ and $|k|\geq N$ there exists a unique eigenvalue of
$L_{t}(Q)$ lying in $\Delta_{k,l}.$ We denote this unique eigenvalue by
$\lambda_{k,j}(t).$ It follows from $(b)$ that the eigenvalues $\lambda
_{k,j}(t)$ for $|k|\geq N$ are simple and the formulas (36), (37) hold.

It remains to prove that the root functions $\Psi_{k,j,p}(x)$ of $L_{t}\left(
Q\left(  x\right)  \right)  $ form a \ Riesz basis, where $\Psi_{k,j,0}(x)$%
{}{}$=\Psi_{k,j}(x)$ and $\Psi_{k,j,1}(x),$ $\Psi_{k,j,2}(x),...,$ are
associated functions belonging to the eigenfunction $\Psi_{k,j}(x).$ Since the
eigenvectors $v_{1},v_{2},...,v_{m}$ of the matrix $C$ form a basis in
$\mathbb{C}^{m}$ and $\mid v_{j}\mid=1$ for $j=1,2,...,m$ the asymptotic
formula (37) yields
\begin{equation}
\sum_{k=1}^{N}(\sum_{p=0,1,...}\sum_{j=1}^{m}\left\vert (f,\Psi_{k,j,p}%
(x))\right\vert ^{2})+\sum_{j=1}^{m}\sum\limits_{k=N+1}^{\infty}\left\vert
(f,\Psi_{k,j}(x))\right\vert ^{2}<\infty.
\end{equation}
for every $f(x)\in L_{2}^{m}(0,1).$ Let $\{\chi_{k,j,p}\}$ be the biorthogonal
to $\{\Psi_{k,j,p}\}$ system of eigenfunctions and associated functions of
$L_{t}^{\ast}(Q).$ Clearly
\[
\chi_{k,j}(x)=\frac{\Psi_{k,j}^{\ast}(x)}{(\Psi_{k,j}(x),\Psi_{k,j}^{\ast
}(x))}%
\]
for $|k|\geq N.$ Hence using (37), (40), we get $(\Psi_{k,j}(x),\Psi
_{k,j}^{\ast}(x))=(v_{j},v_{j}^{\ast})\neq0$ and%
\begin{equation}
\sum_{k=1}^{N}(\sum_{p=0,1,...}\sum_{j=1}^{m}\left\vert (f,\chi_{k,j,p}%
(x))\right\vert ^{2})+\sum_{j=1}^{m}\sum\limits_{k=N+1}^{\infty}\left\vert
(f,\chi_{k,j}(x))\right\vert ^{2}<\infty.
\end{equation}
Since $\{\Psi_{k,j,p}(x)\}$ and the biorthogonal system $\{\chi_{k,j,p}%
\}$\ are total (see [5], where it is proved that the root functions of the
regular boundary value problem in the space of the vector functions form a
Riesz basis with brackets), by the well-known theorem of Bari (see [2], chap.
6), \ the inequalities (41) and (42) imply that the system of the
eigenfunctions and associated functions of $L_{t}(Q)$ forms a Riesz basis in
$L_{2}^{m}(0,1)$
\end{proof}

\end{document}